\documentclass{article}
\usepackage{amsfonts}
\usepackage{hyperref}

\newtheorem{theorem}{Theorem}

\newtheorem{proposition}[theorem]{Proposition}
\newtheorem{remark}[theorem]{Remark}

\begin{document}

\title{Lagrange-Hamilton geometry applied to a Lotka-Volterra dynamical
system}
\author{Ana-Maria Boldeanu and Mircea Neagu}
\date{}
\maketitle

\begin{abstract}
The aim of this paper is to develop, via the least squares variational
method, the Lagrange-Hamilton geometry (in the sense of nonlinear
connections, d-torsions and Lagrangian Yang-Mills electromagnetic-like
energy) produced by a Lotka-Volterra dynamical system, a simple model of the
population dynamics of species competing for some common resource. From a
geometrical point of view, the Jacobi stability of this system is discussed.
\end{abstract}



\section{The Lotka-Volterra dynamical system}

A community of three mutually competing species is modeled by the
Lotka-Volterra dynamical system \cite{Munteanu}:

\begin{equation}
\left\{ 
\begin{array}{l}
\displaystyle\dot{x}_{1}=x_{1}\left(
b_{1}-a_{11}x_{1}-a_{12}x_{2}-a_{13}x_{3}\right) \medskip  \\ 
\displaystyle\dot{x}_{2}=x_{2}\left(
b_{2}-a_{21}x_{1}-a_{22}x_{2}-a_{23}x_{3}\right) \medskip  \\ 
\displaystyle\dot{x}_{3}=x_{3}\left(
b_{3}-a_{31}x_{1}-a_{32}x_{2}-a_{33}x_{3}\right) ,%
\end{array}%
\right.   \label{PR}
\end{equation}%
where $x_{i}(t)$, $i=\overline{1,3}$, is the population size of the $i$-th
species at time $t$, $\dot{x}_{i}$ denote $dx_{i}/dt$, and $a_{ij}$, $b_{i}$%
, $i,j=\overline{1,3}$, are some strictly positive real numbers.

\begin{remark}
The systems of differential equations that represent the competition between
species, that is the so-called Lotka-Volterra systems, were initially
introduced by Lotka \cite{Lotka} and Volterra \cite{Volterra1}, \cite%
{Volterra2}.
\end{remark}

\section{From Lotka-Volterra dynamical system to Lagrange geometry and Jacobi stability}

We consider that the Lotka-Volterra dynamical system (\ref{PR}) is working
on the $3$-dimensional manifold $M=\mathbb{R}^{3},$ whose coordinates are $%
\left( x^{1},x^{2},x^{3}\right) $. Moreover, the corresponding tangent $TM$
and cotangent $T^{\ast }M$ bundles have the coordinates\ $\left(
x^{i},y^{i}\right) _{i=\overline{1,3}}$, respectively $\left(
x^{i},p_{i}\right) _{i=\overline{1,3}}$.

Consequently, if we take the vector field $X=\left( X^{i}\left(
x^{1},x^{2},x^{3}\right) \right) _{i=\overline{1,3}}$ on the manifold $M=%
\mathbb{R}^{3}$, which is expressed by%
\begin{equation}
\begin{array}{l}
\displaystyle X^{1}\left( x^{1},x^{2},x^{3}\right) =x_{1}\left(
b_{1}-a_{11}x_{1}-a_{12}x_{2}-a_{13}x_{3}\right) ,\medskip  \\ 
\displaystyle X^{2}\left( x^{1},x^{2},x^{3}\right) =x_{2}\left(
b_{2}-a_{21}x_{1}-a_{22}x_{2}-a_{23}x_{3}\right) ,\medskip  \\ 
\displaystyle X^{3}\left( x^{1},x^{2},x^{3}\right) =x_{3}\left(
b_{3}-a_{31}x_{1}-a_{32}x_{2}-a_{33}x_{3}\right) ,%
\end{array}
\label{field}
\end{equation}%
we can regard the Lotka-Volterra dynamical system (\ref{PR}) as the
dynamical system%
\begin{equation}
\frac{dx^{i}}{dt}=X^{i}(x(t)),\quad i=\overline{1,3}.  \label{DS}
\end{equation}

It is obvious now that the solutions of class $C^{2}$ of the dynamical
system (\ref{DS}) are the global minimum points for the \textit{least
squares Lagrangian} $L:TM\rightarrow \mathbb{R}$, which is defined by%
\footnote{%
The Latin indices $i,j,k,...$ run from $1$ to $3.$ Moreover, the Einstein
convention of summation is adopted all over this work.}%
\begin{equation}
\begin{array}{c}
L(x,y)=\delta _{ij}\left( y^{i}-X^{i}(x)\right) \left( y^{j}-X^{j}(x)\right)
\Leftrightarrow \medskip  \\ 
L(x,y)=\left( y^{1}-X^{1}(x)\right) ^{2}+\left( y^{2}-X^{2}(x)\right)
^{2}+\left( y^{3}-X^{3}(x)\right) ^{2}.%
\end{array}
\label{LSL}
\end{equation}%
The Euler-Lagrange equations of the least squares Lagrangian (\ref{LSL}) are
described by%
\[
\frac{\partial L}{\partial x^{k}}-\frac{d}{dt}\left( \frac{\partial L}{%
\partial y^{k}}\right) =0,\quad y^{k}=\frac{dx^{k}}{dt},\quad k=\overline{1,3%
},
\]%
and they can be rewritten in the geometrical form:%
\begin{equation}
\frac{d^{2}x^{k}}{dt^{2}}+2G^{k}(x,y)=0,\quad k=\overline{1,3},  \label{E-L}
\end{equation}%
where ($k=\overline{1,3})$%
\[
G^{k}=\frac{1}{2}\left( \frac{\partial ^{2}L}{\partial x^{j}\partial y^{k}}%
y^{j}-\frac{\partial L}{\partial x^{k}}\right) =-\frac{1}{2}\left[ \left( 
\frac{\partial X^{k}}{\partial x^{j}}-\frac{\partial X^{j}}{\partial x^{k}}%
\right) y^{j}+\frac{\partial X^{j}}{\partial x^{k}}X^{j}\right] 
\]%
is endowed with the geometrical meaning of \textit{semispray} of $L$.

In the sequel, following the geometrical ideas from the works Miron and
Anastasiei \cite{Mir-An}, Udri\c{s}te and Neagu \cite{Udr}, \cite{Nea-Udr},
and Balan and Neagu \cite{Bal-Nea}, we construct an entire natural
collection of nonzero Lagrangian geometrical objects (such as nonlinear
connection, d-torsions and Yang-Mills electromagnetic-like energy) that
characterize the Lotka-Volterra dynamical system (\ref{DS}) and, implicitly,
the system (\ref{PR}).

For instance, let us compute the Jacobian matrix $J=\left( \frac{\partial
X^{i}}{\partial x^{j}}\right) _{i,j=\overline{1,3}}$ of the vector field $X$:%
\[
J=\left( 
\begin{array}{ccc}
J_{11} & a_{12}x_{1} & a_{13}x_{1} \\ 
a_{21}x_{2} & J_{22} & a_{23}x_{2} \\ 
a_{31}x_{3} & a_{32}x_{3} & J_{33}%
\end{array}%
\right) ,
\]%
where 
\[
J_{11}=b_{1}-2a_{11}x_{1}-a_{12}x_{2}-a_{13}x_{3},\quad
J_{22}=b_{2}-a_{21}x_{1}-2a_{22}x_{2}-a_{23}x_{3},
\]%
\[
J_{33}=b_{3}-a_{31}x_{1}-a_{32}x_{2}-2a_{33}x_{3}.
\]

It follows that the entries of the canonical nonlinear connection%
\[
N=\left( N_{j}^{i}\right) _{i,j=\overline{1,3}}
\]%
on the manifold $\mathbb{R}^{3}$, which are given by the formulas $%
N_{j}^{i}=\partial G^{i}/\partial y^{j}$ (see \cite{Mir-An}), imply the
following matrix relation (see \cite{Bal-Nea})%
\[
N=-\displaystyle{\frac{1}{2}}\left[ J-J^{t}\right] .
\]%
Consequently, the following result holds good:

\begin{proposition}
The \textbf{Lagrangian canonical nonlinear connection} on the tangent bundle 
$TM$, produced by the Lotka-Volterra dynamical system (\ref{PR}), has the
form:%
\[
N=\left( N_{j}^{i}\right) _{i,j=\overline{1,3}}=\left( 
\begin{array}{ccc}
0 & N_{2}^{1} & N_{3}^{1} \\ 
-N_{2}^{1} & 0 & N_{3}^{2} \\ 
-N_{3}^{1} & -N_{3}^{2} & 0%
\end{array}%
\right) ,
\]%
where%
\[
N_{2}^{1}=-\frac{1}{2}\left( a_{12}x_{1}-a_{21}x_{2}\right) ,\quad
N_{3}^{1}=-\frac{1}{2}\left( a_{13}x_{1}-a_{31}x_{3}\right) ,
\]%
\[
N_{3}^{2}=-\frac{1}{2}\left( a_{23}x_{2}-a_{32}x_{3}\right) .
\]
\end{proposition}

Further, the canonical Cartan linear connection on $TM$ produced by the
least squares Lagrangian (\ref{LSL}) has all adapted components equal to
zero, and the general formulas which give its Lagrangian d-torsions are (see 
\cite{Mir-An})%
\[
R_{k}=\left( R_{jk}^{i}:={\frac{\delta N_{j}^{i}}{\delta x^{k}}}-\frac{%
\delta N_{k}^{i}}{\delta x^{j}}\right) _{i,j=\overline{1,3}},
\]%
where%
\[
{\frac{\delta }{\delta x^{k}}={\frac{\partial }{\partial x^{k}}}-N_{k}^{r}{%
\frac{\partial }{\partial y^{r}}.}}
\]%
By direct computations, it follows that we have (see \cite{Bal-Nea})%
\[
R_{k}=\frac{\partial N}{\partial x^{k}},\quad \forall \;k=\overline{1,3}.
\]%
As a consequence, we find

\begin{proposition}
The Lagrangian canonical Cartan linear connection, produced by the
Lotka-Volterra dynamical system (\ref{PR}), is characterized by the
following $d-$\textbf{torsion skew-symmetric matrices}:%
\[
R_{1}=-\frac{1}{2}\left( 
\begin{array}{ccc}
0 & a_{12} & a_{13} \\ 
-a_{12} & 0 & 0 \\ 
-a_{13} & 0 & 0%
\end{array}%
\right) ,\quad R_{2}=-\frac{1}{2}\left( 
\begin{array}{ccc}
0 & -a_{21} & 0 \\ 
a_{21} & 0 & a_{23} \\ 
0 & -a_{23} & 0%
\end{array}%
\right) , 
\]%
\[
R_{3}=-\frac{1}{2}\left( 
\begin{array}{ccc}
0 & 0 & -a_{31} \\ 
0 & 0 & -a_{32} \\ 
a_{31} & a_{32} & 0%
\end{array}%
\right) . 
\]
\end{proposition}

\begin{proposition}
The\textbf{\ Lagrangian Yang-Mills e\-lec\-tro\-mag\-ne\-tic-like energy}%
\textit{, produced by }the Lotka-Volterra dynamical system (\ref{PR}), is
given by%
\[
\mathcal{EYM}(x)={\frac{1}{2}}\cdot Trace\left[ F\cdot F^{t}\right] ,
\]%
where the electromagnetic-like matrix is $F=-N.$ For more details about
these geometric-physical ideas, see the books \cite{Mir-An} and \cite%
{Bal-Nea}.
\end{proposition}

In the Kosambi-Cartan-Chern (KCC) geometrical theory, it is known that the 
\textit{matrix of deviation curvature} is defined by the formula%
\[
P=\left( P_{j}^{i}\right) _{i,j=\overline{1,3}}=\frac{\partial N}{\partial
x^{k}}y^{k}+\mathcal{E},
\]%
where, if%
\[
\mathcal{E}^{i}=2G^{i}-N_{j}^{i}y^{j}=-\frac{1}{2}\left( \frac{\partial X^{i}%
}{\partial x^{j}}-\frac{\partial X^{j}}{\partial x^{i}}\right) y^{j}-\frac{%
\partial X^{j}}{\partial x^{i}}X^{j}
\]%
is the\textit{\ first invariant of the semispray} of the Lagrangian (\ref%
{LSL}), then we have $\mathcal{E=}\left( \delta \mathcal{E}^{i}/\delta
x^{j}\right) _{i,j=\overline{1,3}}.$ It is clear that, by direct
computations, the entries of the matrix $\mathcal{E}$ are given by%
\[
\frac{\delta \mathcal{E}^{i}}{\delta x^{j}}=-\frac{1}{2}\left( \frac{%
\partial ^{2}X^{i}}{\partial x^{j}\partial x^{k}}-\frac{\partial ^{2}X^{k}}{%
\partial x^{i}\partial x^{j}}\right) y^{k}-\frac{\partial ^{2}X^{k}}{%
\partial x^{i}\partial x^{j}}X^{k}-\frac{\partial X^{k}}{\partial x^{i}}%
\frac{\partial X^{k}}{\partial x^{j}}-
\]%
\[
-\frac{1}{4}\left( \frac{\partial X^{k}}{\partial x^{j}}-\frac{\partial X^{j}%
}{\partial x^{k}}\right) \left( \frac{\partial X^{i}}{\partial x^{k}}-\frac{%
\partial X^{k}}{\partial x^{i}}\right) .
\]

\begin{remark}
Here above, the Einstein convention of summation was used. The indices took
values from $1$ to $3.$
\end{remark}

In conclusion, following the geometrical ideas from the works B\"{o}hmer et
al. \cite{Bohmer}, Buc\u{a}taru-Miron \cite{Buc-Mir} and Neagu-Ovsiyuk \cite%
{Nea-Ovs}, the behavior of the neighboring solutions of the Euler-Lagrange
equations (\ref{E-L}) is \textit{Jacobi stable} if and only if the real
parts of the eigenvalues of the deviation tensor $P$ are strictly negative
everywhere, and \textit{Jacobi unstable}, otherwise. It is also known that
the Jacobi stability or instability has the geometrical meaning that the
trajectories of the Euler-Lagrange equations (\ref{E-L}) are bunching
together or are dispersing (or even are chaotic).

\section{From Lotka-Volterra system to Hamilton geometry}

Let us construct the \textit{least squares Hamiltonian }$H:T^{\ast
}M\rightarrow \mathbb{R}$ associated with the Lagrangian (\ref{LSL}), which
is expressed by%
\begin{equation}
\begin{array}{c}
H(x,p)=\displaystyle{\frac{\delta ^{ij}}{4}}p_{i}p_{j}+X^{k}(x)p_{k}%
\Leftrightarrow \medskip  \\ 
H(x,p)=\displaystyle{\frac{1}{4}}\left( p_{1}^{2}+p_{2}^{2}+p_{3}^{2}\right)
+X^{1}(x)p_{1}+X^{2}(x)p_{2}+X^{3}(x)p_{3},%
\end{array}
\label{LSH}
\end{equation}%
where $p_{r}=\partial L/\partial y^{r}$ and $H=p_{r}y^{r}-L$. It follows
that we can develop again a natural and distinct collection of nonzero
Hamiltonian geometrical objects (such as nonlinear connection and
d-torsions), which also characterize the Lotka-Volterra dynamical system (%
\ref{PR}). For all details about Hamilton geometry on cotangent bundles and
the Hamiltonian least squares variational method for dynamical systems, see
the works: Miron et al. \cite{Miron-Hr-Shim-Sab} and Neagu-Oan\u{a} \cite%
{Nea-Oana}.

The Hamiltonian nonlinear connection on the cotangent bundle $T^{\ast }M$,
has the components (see \cite{Miron-Hr-Shim-Sab})%
\[
N_{ij}=\frac{\partial ^{2}H}{\partial x^{j}\partial p_{i}}+\frac{\partial
^{2}H}{\partial x^{i}\partial p_{j}}.
\]%
By direct computations, we get $\mathbf{N}=J+$ $J^{t}$. For more details,
see the monograph \cite{Nea-Oana}. Consequently, we infer

\begin{proposition}
The \textbf{Hamiltonian canonical nonlinear connection} on the cotangent
bundle $T^{\ast }M$, produced by the Lotka-Volterra dynamical system (\ref%
{PR}), is described by the matrix%
\[
\mathbf{N}=\left( N_{ij}\right) _{i,j=\overline{1,6}}=\left( 
\begin{array}{ccc}
N_{11} & N_{12} & N_{13} \\ 
N_{12} & N_{22} & N_{23} \\ 
N_{13} & N_{23} & N_{33}%
\end{array}%
\right) ,
\]
\end{proposition}

where 
\[
N_{11}=2\left( b_{1}-2a_{11}x_{1}-a_{12}x_{2}-a_{13}x_{3}\right) , 
\]%
\[
N_{12}=a_{12}x_{1}+a_{21}x_{2},\quad 
\]%
\[
N_{13}=a_{13}x_{1}+a_{31}x_{3}, 
\]%
\[
N_{22}=2\left( b_{2}-a_{21}x_{1}-2a_{22}x_{2}-a_{23}x_{3}\right) , 
\]%
\[
N_{23}=a_{23}x_{2}+a_{32}x_{3}, 
\]%
\[
N_{33}=2\left( b_{3}-a_{31}x_{1}-a_{32}x_{2}-2a_{33}x_{3}\right) . 
\]

Further, the canonical Cartan linear connection on $T^{\ast }M$ produced by
the least squares Hamiltonian (\ref{LSH}) has all adapted components equal
to zero. Moreover, the general formulas which give the Hamiltonian
d-torsions are (see \cite{Miron-Hr-Shim-Sab})%
\[
\mathbf{R}_{k}=\left( R_{kij}:={\frac{\delta N_{ki}}{\delta x^{j}}}-{\frac{%
\delta N_{kj}}{\delta x^{i}}}\right) _{i,j=\overline{1,3}},
\]%
where%
\[
{\frac{\delta }{\delta x^{j}}={\frac{\partial }{\partial x^{j}}}-N_{rj}{%
\frac{\partial }{\partial p_{r}}.}}
\]%
It follows that we find (for all details, see also \cite{Nea-Oana})

\begin{proposition}
The Hamiltonian canonical Cartan linear connection, produced by the
Lotka-Volterra dynamical system (\ref{PR}), is characterized by the
following $d-$\textbf{torsion skew-symmetric matrices}:%
\[
\mathbf{R}_{k}=\frac{\partial }{\partial x^{k}}\left[ J-J^{t}\right]
=-2R_{k},\quad \forall \;k=\overline{1,3}. 
\]
\end{proposition}

\section{Conclusion}

In the context of this paper, from our new
geometric-physical approach, the surfaces of constant level of the
Lagrangian Yang-Mills electromagnetic-like energy produced by the
Lotka-Volterra dynamical system (\ref{PR}) could have important connotations for the phenomenom taken in study. Thus, the surface of constant
level $\Sigma _{\rho }$ is obviously described by the equation%
\[
\left( a_{12}x_{1}-a_{21}x_{2}\right) ^{2}+\left(
a_{13}x_{1}-a_{31}x_{3}\right) ^{2}+\left( a_{23}x_{2}-a_{32}x_{3}\right)
^{2}=4\rho \geq 0.
\]%
The question is: There is a meaning, related to the Lotka-Volterra competition between species, for the shape of this quadric surface? 

For such a reason, the computer shown graphics of this surface of constant level, together with its differential geometry (in the sense of curvatures, asymptotic lines, curvature lines or geodesics) represents our future work in progress, and it is in our attention.

\noindent \textbf{Open problem. }What is the real meaning in the Lotka-Volterra competition between species for our Lagrange-Hamilton geometrical objects constructed in this paper?

\noindent \textsc{Ana-Maria Boldeanu and Mircea Neagu}\\[1mm]
Transilvania University of Bra\c{s}ov\newline
Department of Mathematics and Computer Science\newline
Blvd. Iuliu Maniu, No. 50, 500091 Bra\c{s}ov, Romania.

\noindent E-mails: ana.boldeanu@unitbv.ro, mircea.neagu@unitbv.ro

\end{document}